\documentclass{mefdreport}

\usepackage{multicol}
\usepackage{multirow}
\usepackage{fancybox}
\usepackage{index}
\usepackage{varioref}
\usepackage{psfrag}
\usepackage{epsfig}
\usepackage{boxedminipage}
\usepackage{graphicx}
\usepackage{rotating}
\usepackage{amsmath}
\usepackage{amssymb}
\usepackage{bm}
\usepackage{latexsym}
\usepackage{alltt}
\usepackage[small,bf]{caption}
\usepackage{url}
\usepackage{citesort}
\usepackage{crop}
\usepackage{array}
\usepackage{subfigure}
\usepackage{dcolumn}
\usepackage{mathrsfs}

\newcommand{\EQ}[1]{(\ref{eq:#1})}

\newcommand{\FIG}[1]{\ref{fig:#1}}
\newcommand{\SEC}[1]{\ref{sec:#1}}

\newcommand{\grad}{\nabla}
\newcommand{\laplace}{\Delta}
\renewcommand{\div}{\nabla\cdot{}}

\newcommand{\dd}{\mathrm{d}}

\newcommand{\IR}{\mathbb{R}}

\newcommand{\optype}[1]{\mathrm{#1}}

\newcommand{\meas}[1]{\optype{meas}}

\newcommand{\domF}{\smash[tb]{\Omega^{\mathrm{f}}_t}}
\newcommand{\domFi}{\smash[tb]{\Omega^{\mathrm{f}}_0}}
\newcommand{\domFref}{\smash[tb]{\hat{\Omega}{}^{\mathrm{f}}}}
\newcommand{\domS}{\smash[tb]{\Omega^{\mathrm{s}}_t}}
\newcommand{\domSi}{\smash[tb]{\Omega^{\mathrm{s}}_0}}
\newcommand{\domSref}{\smash[tb]{\hat{\Omega}{}^{\mathrm{s}}}}
\newcommand{\Gammaref}{\smash[tb]{\hat{\Gamma}}}

\begin{document}

\title{Elasto-capillarity Simulations based on the Navier--Stokes--Cahn--Hilliard Equations}
\author{E.H. van Brummelen, M. Shokrpour-Roudbari and G.J. van Zwieten}
\thanks{%
Eindhoven University of Technology, Department of Mechanical Engineering\\
  P.O. Box 513, 5600 MB Eindhoven, The Netherlands}
\maketitle

\abstract{
We consider a computational model for complex-fluid-solid interaction
based on a diffuse-interface model for the complex fluid and a hyperelastic-material model for the solid. The diffuse-interface complex-fluid model is described by the 
incompressible Navier--Stokes--Cahn--Hilliard equations with preferential-wetting 
boundary conditions at the fluid-solid interface. The corresponding fluid traction
on the interface includes a capillary-stress contribution,
and the dynamic interface condition comprises the traction exerted by the non-uniform fluid-solid surface tension. We present a weak formulation of the aggregated complex-fluid-solid-interaction problem, based on an Arbitrary-Lagrangian-Eulerian formulation
of the Navier--Stokes--Cahn--Hilliard equations and a proper reformulation of the 
complex-fluid traction and the fluid-solid surface tension. To validate the presented
complex-fluid-solid-interaction model, we present numerical results and
conduct a comparison to experimental data for a droplet on a soft substrate.
}

\section{Introduction}
\label{sec:intro}

Complex fluids are fluids that consist of multiple constituents, e.g. of multiple phases of the same fluid (gas, liquid or solid) or of multiple distinct species (e.g. water and air). The interaction of such complex fluids with elastic solids leads to multitudinous intricate physical phenomena. Examples are {\em durotaxis\/}, viz., seemingly spontaneous migration of liquid droplets on solid substrates with an elasticity gradient~\cite{Style30072013}, or {\em capillary origami}, viz., large-scale solid deformations induced by capillary forces~\cite{PhysRevLett.98.156103}. Complex-Fluid-Solid Interaction (CFSI) is moreover of fundamental technological importance in a wide variety of applications, such as inkjet printing and additive manufacturing.

Despite significant progress in models and computational techniques for 
the interaction of solids and classical fluids (see~\cite{Bazilevs:2013zm,Tezduyar:2010ys} for an overview), and for complex fluids separately 
(see, e.g.,~\cite{Lowengrub:1998uq,Jacqmin:1999fk,Hohenberg:1977hh,Abels:2012vn,Guo:2014hb,Aland:2012fk}), complex-fluid-solid interaction has remained essentially unexplored. A notable exception is the computational 
CFSI model based on the Navier--Stokes--Korteweg equations in~\cite{Bueno:2014qv}.

In this contribution we consider a computational model for complex-fluid-solid interaction
based on a diffuse-interface complex-fluid model and a hyperelastic solid model with a Saint Venant--Kirchhoff stored-energy functional. The diffuse-interface complex-fluid model is described by the incompressible Navier--Stokes--Cahn--Hilliard (NSCH) equations.
The interaction of the complex fluid with the solid substrate is represented
by dynamic and kinematic interface conditions and a preferential-wetting 
boundary condition. The traction exerted by the complex fluid on the fluid-solid 
interface comprises a non-standard capillary-stress contribution, in addition 
to the standard pressure and viscous-stress components. The dynamic condition
imposes equilibrium of this complex-fluid traction, the traction exerted by the
hyperelastic solid and the traction due to the non-uniform fluid-solid surface tension.
We present a weak formulation of the aggregated complex-fluid-solid-interaction problem, based on an Arbitrary-Lagrangian-Eulerian (ALE) formulation
of the NSCH system and a suitable weak representation of the complex-fluid 
traction and the non-uniform fluid-solid surface tension. 

To evaluate the capability of the considered complex-fluid-solid-interaction
model to describe elasto-capillary phenomena, we consider numerical experiments
for a test case pertaining to a droplet on a soft substrate, and we present
a comparison to experimental data from~\cite{Style:2013ay}.

The remainder of this contribution is organized as follows. Section~\SEC{probstat}
presents a specification and discussion of the considered complex-fluid-solid-interaction problem. In Section~\SEC{weakform} we treat the weak formulation of the aggregated 
fluid-solid-interaction problem. Section~\SEC{numexp} is concerned with 
numerical experiments and results. Concluding remarks are presented in Section~\SEC{concl}.

\section{Problem Statement}
\label{sec:probstat}
To accommodate the complex-fluid-solid system, we consider a time interval $(0,T)\subseteq\IR_{>0}$ and two simply-connected time-dependent open subsets 
$\smash[tb]{\domF}\subset\IR^d$ ($d=2,3$) and $\smash[tb]{\domS}\subset\IR^d$,
which hold the complex-fluid and solid, respectively.
The fluid-solid interface corresponds to~$\smash[tb]{\Gamma_t:=\partial\domF\cap\partial\domS}\neq\emptyset$. We assume that the time-dependent configuration $\smash[tb]{\Omega_t:=\mathrm{int}(\mathrm{cl}\,\domF\cup\mathrm{cl}\,\domS)}$
is the image of a time-dependent transformation $\hat{\chi}$ 
acting on a fixed reference domain
$\smash[tb]{\hat{\Omega}:=\mathrm{int}(\mathrm{cl}\,\domFref\cup\mathrm{cl}\,\domSref)}$ such that $\domF=\hat{\chi}\domFref$ and $\domS=\hat{\chi}\domSref$. The restrictions of 
${\hat{\chi}}$ to $\domFref$ and $\domSref$ are denoted by ${\hat{\chi}}^{\mathrm{f}}$ and~${\hat{\chi}}^{\mathrm{s}}$, respectively.
The reference domains $\domFref:=\domFi$ and $\domSref:=\domSi$ are identified with the initial configurations.

\subsection{Navier--Stokes--Cahn--Hilliard Complex-Fluid Model}
We consider a complex fluid composed of two immiscible incompressible constituents, separated by a thin diffuse interface. The behavior of the complex fluid is described by the Navier--Stokes--Cahn--Hilliard (NSCH) equations. The two species are identified by an order parameter $\varphi:\domF\to[-1,1]$. Typically, $\varphi$ is selected as either volume fraction~\cite{Hohenberg:1977hh,Abels:2012vn,Jacqmin:1999fk} 
or mass fraction~\cite{Lowengrub:1998uq,Guo:2014hb}, such that $\varphi=1$ (resp. $\varphi=-1$) pertains to a pure species\nobreakdash-1 (resp. species\nobreakdash-2) composition of the fluid, and 
$\varphi\in(-1,1)$ indicates a mixture. Depending on the definition
of the phase indicator~$\varphi$ as mass or volume fraction, and the definition of mixture velocity as mass-averaged or volume-averaged species velocity, various forms of the NSCH equations can be derived. In mass-averaged-velocity formulations, the mixture is generally quasi-incompressible. In volume-averaged-velocity formulations, the mixture is incompressible.
We select~$\varphi$ as volume fraction and consider a volume-averaged-velocity formulation. The behavior of the complex fluid is described by~\cite{Aland:2012fk,Jacqmin:1999fk}:
\begin{equation}
\label{eq:NSCH}
\left.
\begin{aligned}
\partial_t(\rho{u})+\div(\rho{}u\otimes{}u)+\grad{}p-\div\tau+
\tilde{\sigma}\epsilon\div(\grad\varphi\otimes\grad\varphi)&=0
\\
\div{}u&=0
\\
\partial_t{\varphi}+\div(\varphi{}u)-\gamma\laplace\mu&=0
\\
\mu+\tilde{\sigma}\epsilon\laplace\varphi-\tilde{\sigma}\epsilon^{-1}W'(\varphi)&=0
\end{aligned}
\right\}\text{ in }\domF
\end{equation}
with $\rho:=\rho(\varphi)=\rho_1(1+\varphi)/2+\rho_2(1-\varphi)/2$ as mixture density, 
$u:\domF\to\IR^d$~as volume-averaged mixture velocity, 
$p=p:\domF\to\IR$ as pressure,
$\tau=\nu\grad^{s}u$ as viscous-stress tensor and
$\mu:\domF\to\IR$ as chemical potential. The mixture viscosity is defined as $\nu:=\nu(\varphi)=\nu_1(1+\varphi)/2+\nu_2(1-\varphi)/2$. The 
parameter~$\tilde{\sigma}$ is related to the fluid-fluid surface tension $\sigma$
by $2\sqrt{2}\,\tilde{\sigma}=3\,\sigma$, and $\gamma>0$ designates mobility. 
The energy density associated with mixing of the constituents is represented by the standard double-well potential $W(\varphi)=\frac{1}{4}(\varphi^2-1)^2$. The parameter~$\epsilon>0$ controls the thickness of the diffuse interface between the fluid constituents.

Suitable initial conditions for~\EQ{NSCH} are provided by
a specification of the initial phase distribution and the initial velocity, according to $\varphi(0,\cdot)=\varphi_0$ and $u(0,\cdot)=u_0$, with $\varphi_0:\domFi\to[-1,1]$ and $u_0:\domFi\to\IR^d$ exogenous data. Equations~(\ref{eq:NSCH}$_{1}$),
(\ref{eq:NSCH}$_3$) and~(\ref{eq:NSCH}$_4$) are 
typically furnished with Dirichlet or Neumann boundary conditions:
\enlargethispage{\baselineskip}
\begin{equation}
\label{eq:NSCHBC}
\begin{aligned}
  u&=g{}^u_{\mathrm{D}}&&\text{on }\Gamma^u_{\mathrm{D}}
\\
\varphi&=g^{\varphi}_{\mathrm{D}}&&\text{on }\Gamma^{\varphi}_{\mathrm{D}}
  \\
\mu&=g{}^{\mu}_{\mathrm{D}}&&\text{on }\Gamma^{\mu}_{\mathrm{D}}
\end{aligned}
\qquad\qquad
\begin{aligned}
-pn+\tau{}n-\tilde{\sigma}\epsilon\partial_n\varphi\grad\varphi&=g{}^u_{\mathrm{N}}
&&\text{on }\Gamma^u_{\mathrm{N}}
\\
-\tilde{\sigma}\epsilon\partial_n\phi&=g{}^{\varphi}_{\mathrm{N}}
&&\text{on }\Gamma^{\varphi}_{\mathrm{N}}
\\
\gamma\partial_n\mu&=g{}^{\mu}_{\mathrm{N}}
&&\text{on }\Gamma^{\mu}_{\mathrm{N}}
\end{aligned}
\end{equation}
with $n$ the exterior unit normal vector to~$\partial\domF$.
The right-hand sides in~\EQ{NSCHBC} correspond to exogenous data.
It generally holds that $\smash[tb]{\Gamma_{\mathrm{D}}^{(\cdot)}\cap\Gamma_{\mathrm{N}}^{(\cdot)}}=\emptyset$. The Neumann condition in~(\ref{eq:NSCHBC}$_1$) provides a specification of the fluid traction on the 
boundary~$\smash[tb]{\Gamma_{\mathrm{N}}^u}$. 
If $\smash[tb]{\Gamma_{\mathrm{N}}^{\mu}}$ corresponds to a material boundary,
homogeneous data $\smash[tb]{g_{\mathrm{N}}^{\mu}}=0$ provide
phase conservation at the boundary.
Indeed, from~(\ref{eq:NSCH}$_3$), the Neumann condition 
in~(\ref{eq:NSCHBC}$_3$) and the Reynolds transport theorem it follows that:
\begin{equation}
\label{eq:phaseconservation}
\frac{d}{dt}\int_{\domF}\varphi
=
\int_{\partial\domF\setminus\Gamma_{\mathrm{N}}^{\mu}}
\big(\gamma\partial_n\mu-\varphi(u_n-w_n)\big)
-
\int_{\Gamma_{\mathrm{N}}^{\mu}}
\varphi(u_n-w_n)
+\int_{\Gamma_{\mathrm{N}}^{\mu}}g_{\mathrm{N}}^{\mu}
\end{equation}
with $u_n$ and $w_n$ the normal velocities of the fluid and of the boundary, respectively. 
Material boundaries satisfy $u_n=w_n$ and, accordingly, the penultimate term in~\EQ{phaseconservation} vanishes. Therefore, 
the contribution of~$\smash[tb]{\Gamma_{\mathrm{N}}^{\mu}}$
to production of~$\varphi$ vanishes if $\smash[tb]{g{}_{\mathrm{N}}^{\mu}}=0$.
An important alternative to~(\ref{eq:NSCHBC}$_2$), is the nonlinear Robin-type condition
\begin{equation}
\label{eq:wetting}
\tilde{\sigma}\epsilon\partial_n\varphi+\sigma_{\mathrm{W}}'(\varphi)=0
\quad\text{on }\Gamma_{\mathrm{W}}
\end{equation} 
with $\sigma_{\mathrm{W}}(\varphi)=\tfrac{1}{4}(\varphi^3-3\varphi)(\sigma_2-\sigma_1)+\tfrac{1}{2}(\sigma_1+\sigma_2)$ the surface tension of the 
complex-fluid-solid interface, and $\sigma_1>0$ and $\sigma_2>0$ the fluid-solid surface tensions of species~1 and~2, respectively; see~\cite{Jacqmin:2000kx}. 
Note that $\sigma_{\mathrm{W}}(\cdot)$ provides an interpolation of the pure-species fluid-solid surface tensions, i.e. $\sigma_{\mathrm{W}}(1)=\sigma_1$ and $\sigma_{\mathrm{W}}(-1)=\sigma_2$. Equation~\EQ{wetting} describes preferential wetting of~$\smash[tb]{\Gamma_{\mathrm{W}}}$ by the two fluid components. In particular, 
the angle $\theta_{\mathrm{s}}=\arccos((\sigma_2-\sigma_1)/\sigma)$ corresponds to the static contact angle between the diffuse interface and $\Gamma_{\mathrm{W}}$ (interior to fluid~1). Interaction of the complex fluid~\EQ{NSCH} 
with a solid substrate is modeled by Dirichlet condition~(\ref{eq:NSCHBC}$_{1}$), 
(homogeneous) Neumann condition~(\ref{eq:NSCHBC}$_{3}$), 
and preferential-wetting condition~\EQ{wetting}.

\subsection{Hyperelastic Saint~Venant--Kirchhoff Solid Model}
We consider a hyperelastic solid with 
Saint Venant--Kirchhoff stress-strain relation. 
Denoting the initial density of the solid by $\hat{\rho}:=\rho_0:\domSref\to\IR_{>0}$, 
the solid deformation $\hat{\chi}^{\mathrm{s}}:\domSref\to\domS$ satisfies
the equation of motion:
\begin{equation}
\label{eq:Mom}
\hat{\rho}\partial_{t}^2\hat{\chi}^{\mathrm{s}}-\hat{\nabla}\cdot{}\hat{P}=0\quad\text{ in }\domSref
\end{equation} 
with $\hat{P}$ the first Piola--Kirchhoff stress tensor and $\hat{\grad}\cdot$ the divergence operator in the reference configuration. For hyperelastic materials, $\hat{\grad}\cdot\hat{P}$ is the vector-valued function such that
$-\int_{\domSref}\hat{x}\cdot(\hat{\grad}\cdot{}\hat{P})
=
\mathcal{W}'(\hat{\chi}^{\mathrm{s}};\hat{x})$
for all 
$\hat{x}\in{}C^{\infty}_0(\domSref,\IR^d)$,
with $\mathcal{W}$ the stored-energy functional, $\mathcal{W}'(\hat{\chi}^{\mathrm{s}};\cdot)$ its Fr\'echet derivative at~$\hat{\chi}^{\mathrm{s}}$,
 and $\smash[tb]{C^{\infty}_0(\domSref,\IR^d)}$ 
the class of $\IR^d$-valued smooth functions with compact support in $\domSref$. Denoting by $F:=F(\hat{\chi}^{\mathrm{s}})$ the deformation tensor and by $E:=\tfrac{1}{2}(F^TF-I)$
the Green--Lagrange strain tensor, the Saint Venant--Kirchhoff relation
specifies the strain-energy density associated with~$\hat{\chi}^{\mathrm{s}}$ as
$\tfrac{1}{2}\lambda_{\mathrm{L}}(\mathrm{tr}\,E)^2+\mu_{\mathrm{L}}(\mathrm{tr}\,E^2)$
with $\lambda_{\mathrm{L}}$ and $\mu_{\mathrm{L}}$ the Lam\'e parameters.

The identification of the reference configuration and the
initial configuration yields the initial condition $\hat{\chi}^{\mathrm{s}}_0=\mathrm{Id}$. 
Equation~\EQ{Mom} 
is generally furnished with Dirichlet or Neumann conditions:
\begin{equation}
\label{eq:SolBC}
\hat{\chi}^{\mathrm{s}}=g^{\hat{\chi}}_{\mathrm{D}}\quad\text{on }\hat{\Gamma}_{\mathrm{D}}^{\hat{\chi}}
\qquad
\hat{P}\hat{n}=g^{\hat{\chi}}_{\mathrm{N}}\quad\text{on }\hat{\Gamma}_{\mathrm{N}}^{\hat{\chi}}
\end{equation}
with $\smash[tb]{g^{\hat{\chi}}_{\mathrm{D}}}$ and 
$\smash[tb]{g^{\hat{\chi}}_{\mathrm{N}}}$ deformation and traction data
on 
$\smash[tb]{\hat{\Gamma}_{\mathrm{D}}^{\hat{\chi}}}$ and
$\smash[tb]{\hat{\Gamma}_{\mathrm{N}}^{\hat{\chi}}}$, respectively.


\subsection{Interface Conditions}
The complex fluid~\EQ{NSCH} and solid~\EQ{Mom} are interconnected at the interface by
kinematic and dynamic interface conditions.
The kinematic condition identifies the mixture velocity and the structural 
velocity at the interface. This condition can be interpreted as a Dirichlet
boundary condition for fluid velocity in accordance with~(\ref{eq:NSCHBC}$_1$):
%
\begin{equation}
\label{eq:kin}
\qquad{}u=g^u_{\mathrm{D}}:=\partial_t\hat{\chi}^{\mathrm{s}}\circ{}\hat{\chi}^{-1}\text{ on }
\Gamma_t\subseteq\Gamma_{\mathrm{D}}^u
\end{equation}
Kinematic condition~\EQ{kin} constitutes a partial solid-wall condition for~\EQ{NSCH}.
The condition is complemented by a homogeneous Neumann
condition~(\ref{eq:NSCHBC}$_3$) to impose conservation of phase,
and wettability boundary condition~(\ref{eq:wetting}).

The dynamic condition imposes equilibrium of the fluid and solid tractions
and the traction exerted on the interface by the fluid-solid surface tension.
The traction due to the fluid-solid surface tension is given by 
the Young-Laplace relation for non-uniform surface tension according to
$\Sigma^{\mathrm{fs}}=\sigma_{\mathrm{W}}(\varphi)\kappa{}\,n+\grad_{\Gamma}\sigma_{\mathrm{W}}(\varphi)$, with $\kappa$ as the additive curvature of~$\Gamma_t$
and $\grad_{\Gamma}(\cdot)$ the tangential gradient on~$\Gamma_t$; see, 
e.g.,~\cite{Gouin:2014ai}. We adopt the convention that curvature is negative if the osculating circle in the normal plane is located in
the fluid domain. The complex fluid in~\EQ{NSCH} exerts traction 
$\smash[tb]{\Sigma{}^{\mathrm{f}}}:=pn-\tau{}n+\tilde{\sigma}\epsilon\partial_n\varphi\grad\varphi$ on the interface; cf.~(\ref{eq:NSCHBC}$_1$). Note the 
capillary-stress contribution, $\tilde{\sigma}\epsilon\partial_n\varphi\grad\varphi$, to
the fluid traction.
The traction exerted by the solid~\EQ{Mom} 
is~$\smash[tb]{\hat{\Sigma}{}^{\mathrm{s}}}:=-\hat{P}\hat{n}$ with $\hat{n}$ the exterior unit normal vector to $\partial\domSref$; cf.~\EQ{SolBC}. To account for the fact that fluid traction 
and surface-tension traction are 
expressed in the current configuration and solid traction is expressed in the reference configuration, we consider the dynamic
condition in distributional form:
\begin{equation}
\label{eq:dyn}
\int_{\hat{\Gamma}}\hat{v}\cdot{}\hat{\Sigma}{}^{\mathrm{s}}\,\dd{}\hat{s}
=
-
\int_{\Gamma_t}\big(\hat{v}\circ{}\hat{\chi}^{-1}\big)\cdot\big(\Sigma{}^{\mathrm{f}}+\Sigma{}^{\mathrm{fs}}\big)
\,\dd{}s
\quad\forall{}\hat{v}\in{}C^{\infty}_0(\hat{\Gamma})
\end{equation}
with $\dd\hat{s}$ and $\dd{}s$ the surface measures carried by $\hat{\Gamma}$ and 
$\Gamma_t$, respectively. A precise interpretation of~\EQ{dyn} based on
weak traction evaluation is presented in Section~\SEC{weakfsi}.

\section{Weak Formulation}
\label{sec:weakform}

In this section we present a consistent weak formulation of the fluid-solid-interaction problem in Section~\SEC{probstat}. We first consider a weak Arbitrary-Lagrang\-ian-Eulerian (ALE) formulation of the NSCH system~\EQ{NSCH} in Section~\SEC{NSCHALE}.
Section~\SEC{weakfsi} presents a weak formulation of aggregated FSI problem,
including a weak formulation of the solid subsystem~\EQ{Mom}
and an appropriate weak formulation of the traction exerted by the fluid
on the solid at the interface in conformity with the dynamic condition. 

\subsection{ALE Formulation of NSCH Equations}
\label{sec:NSCHALE}
To accommodate the motion of the fluid domain, we consider a weak
formulation of~\EQ{NSCH} in ALE form. 
The weak formulation is set in the current
configuration. The deformation of the fluid domain,~$\hat{\chi}^{\mathrm{f}}$, 
induces domain velocity $\hat{w}:=\partial_t\hat{\chi}^{\mathrm{f}}$. To derive the ALE formulation, we note
that
\begin{equation}
\label{eq:ALEID}
\int_{\domF}z\partial_t{}\psi
=
\frac{d}{dt}\int_{\domF}z\psi
-\int_{\domF}\div{}\big(\psi{}w{}z)
-\int_{\domF}\psi\partial_tz
=
\frac{d}{dt}\int_{\domF}z\psi
-\int_{\domF}z\div(\psi{}w)
\end{equation}
for all $\smash[tb]{\hat{z}\in{}C^{\infty}(\domFref)}$ 
and $\smash[tb]{\psi\in{}C^{\infty}(\domF)}$,
with $z=\hat{z}\circ{}\hat{\chi}^{-1}$ and $w=\hat{w}\circ{}\hat{\chi}^{-1}$.
The identities in~\EQ{ALEID} follow from the transport theorem 
and~$\partial_t(\hat{z}\circ{}\hat{\chi}^{-1})=-w\cdot\grad{}z$. 
From~\EQ{ALEID} it follows that~(\ref{eq:NSCH}$_1$) and~(\ref{eq:NSCH}$_3$) 
subject to~\EQ{NSCHBC} can be recast into the weak ALE form:
\begin{equation}
\label{eq:NSCHW1}
\begin{aligned}
d_t\langle\rho{}u,v\rangle+
\mathcal{A}_{\mathrm{N}}(u,w,\varphi;v)+\mathcal{B}(p,v)
&=\mathcal{L}_{\mathrm{N}}(u,w,\varphi,p;v)
&\quad&\forall\hat{v}\in{C}^{\infty}(\domFref,\IR^d)
\\
d_t\langle\varphi,z\rangle
+\mathcal{A}_{\mathrm{C}}(u,w,\varphi,\mu;z)
&=\mathcal{L}_{\mathrm{C}}(u,w,\varphi,\mu;z)
&\quad&\forall{}\hat{z}\in{C}^{\infty}(\domFref)
\end{aligned}
\end{equation}
with $\langle{}\rho{}u,v\rangle=\int_{\domF}v\cdot\rho{}u$ and 
$\langle\varphi,z\rangle=\int_{\domF}z\,\varphi$, and
\begin{align}
\mathcal{A}_{\mathrm{N}}(u,w,\varphi;v)
&=
\int_{\domF}\grad{}v:\big(\tau-\rho{}u\otimes(u-w)
-
\tilde{\sigma}\epsilon\grad\varphi\otimes\grad\varphi\big)
\notag
\\
\mathcal{L}{}_{\mathrm{N}}(u,w,\varphi,p;v)
&=
\int_{\Gamma_{\mathrm{N}}^u}v\cdot{}g_{\mathrm{N}}^u
-\int_{\partial\domF\setminus\Gamma_{\mathrm{N}}^u}v\cdot{}\Sigma^{\mathrm{f}}
-\int_{\partial\domF}v\cdot\rho{}u(u_n-w_n)
\notag
\\
\mathcal{A}_{\mathrm{C}}(u,w,\varphi,\mu;z)
&=
\int_{\domF}\grad{}z\cdot\big(\gamma\grad\mu-\varphi(u-w)\big)
\label{eq:NSCHforms1}
\\
\mathcal{L}_{\mathrm{C}}(u,w,\varphi,\mu;z)
&=
\int_{\Gamma_{\mathrm{N}}^{\mu}}z\,g_{\mathrm{N}}^{\mu}
+\int_{\partial\domF\setminus\Gamma_{\mathrm{N}}^{\mu}}z\,\gamma\partial_n\mu
-\int_{\partial\domF}z\,\varphi(u_n-w_n)
\notag
\\
\mathcal{B}(p,v)&=\int_{\domF}-p\,\div{}v
\notag
\end{align}
%
From~(\ref{eq:NSCH}$_2$) and~(\ref{eq:NSCH}$_4$), the Neumann condition
in~(\ref{eq:NSCHBC}$_2$) and the wetting
condition~(\ref{eq:wetting}), we moreover infer:
\begin{equation}
\label{eq:NSCHW2}
\begin{aligned}
\mathcal{B}(q,u)&=0
&\quad&\forall{}q\in{C}^{\infty}(\domF)
\\
\mathcal{A}_{\mathrm{P}}(\varphi;\mu,y)&=\mathcal{L}_{\mathrm{P}}(\varphi,y)
&\quad&\forall{}y\in{C}^{\infty}(\domF)
\end{aligned}
\end{equation}
with 
\begin{equation}
\label{eq:NSCHforms2}
\begin{aligned}
\mathcal{A}_{\mathrm{P}}(\varphi;\mu,y)&=
\int_{\domF}y\big(\mu-\tilde{\sigma}\epsilon^{-1}W'(\varphi)\big)
-
\int_{\domF}\grad{}y\cdot\tilde{\sigma}\epsilon\grad\varphi
-
\int_{\Gamma_{\mathrm{W}}}y{}f'(\varphi)
\\
\mathcal{L}_{\mathrm{P}}(\varphi,y)&=
\int_{\Gamma_{\mathrm{N}}^{\varphi}}y\,g^{\varphi}_{\mathrm{N}}
-\int_{\partial\domF\setminus(\Gamma_{\mathrm{W}}\cup\Gamma_{\mathrm{N}}^{\varphi})}
y\,\tilde{\sigma}\epsilon\partial_n\varphi
\end{aligned}
\end{equation}
It is important to note that the fluid-solid interface satisfies
$\Gamma_t\subseteq\Gamma_{\mathrm{W}}\cap\Gamma_{\mathrm{N}}^{\mu}\cap\Gamma_{\mathrm{D}}^u$.

The configuration of the fluid domain, $\hat{\chi}^{\mathrm{f}}$, can be constructed in 
various manners. We select $\hat{\chi}^{\mathrm{f}}=
h_{\hat{\chi}^{\mathrm{s}}|_{\hat{\Gamma}}}$ as the harmonic extension of the trace of the solid displacement on the interface onto~$\domFref$. Accordingly, it holds that $\hat{w}=\partial_t{}h_{\hat{\chi}^{\mathrm{s}}|_{\hat{\Gamma}}}$.

\subsection{Aggregated Fluid-Solid-Interaction Problem}
\label{sec:weakfsi}
From the equation of motion of the solid in~\EQ{Mom}, we infer the weak formulation:
\begin{equation}
\label{eq:WeakMom}
d_t^2\langle\hat{\rho}\hat{\chi}^{\mathrm{s}},\hat{x}\rangle
+
\mathcal{W}'(\hat{\chi}^{\mathrm{s}};\hat{x})
=
\int_{\hat{\Gamma}_{\mathrm{N}}^{\hat{\chi}}}\hat{x}\cdot{}g_{\mathrm{N}}^{\hat{\chi}}
-\int_{\partial\domSref\setminus(\Gammaref\cup\hat{\Gamma}_{\mathrm{N}}^{\hat{\chi}})}\hat{x}\cdot\hat{\Sigma}{}^s
-
\int_{\Gammaref}\hat{x}\cdot\hat{\Sigma}{}^s
\end{equation}
for all $\hat{x}\in{}C^{\infty}(\domSref,\IR^d)$.
The ultimate term in~\EQ{WeakMom} constitutes the
solid traction on the interface. The dynamic condition
imposes that this term coincides with the right-hand side of~\EQ{dyn}.
Noting that $\Gamma_t\subseteq\partial\domF\setminus\Gamma_{\mathrm{N}}^u$,
equations~(\ref{eq:NSCHW1})--(\ref{eq:NSCHforms1}) convey that the 
fluid-traction contribution can be expressed as:
\begin{equation}
\label{eq:WeakTrac}
-\int_{\Gamma_t}x\cdot\Sigma{}^{\mathrm{f}}
=
\int_{\Gamma_t}\ell_x\cdot\rho{}u(u_n-w_n)
+
d_t\langle\rho{}u,\ell_x\rangle
+
\mathcal{A}_{\mathrm{N}}(u,w,\varphi;\ell_x)
+
\mathcal{B}(p,\ell_x)
\end{equation}
where $\ell_x$ represents an appropriate lifting of~$x$, viz. any suitable 
function $\domF\to\IR^d$ such that $\ell_x|_{\Gamma_t}=x$ 
and $\ell_x$ vanishes on~$\partial\domF\setminus\Gamma_t$. 
The right member of~\EQ{WeakTrac} provides a weak formulation
of the traction functional in the left member of~\EQ{WeakTrac},
in the sense that the identity~\EQ{WeakTrac} holds for all solutions 
of~\EQ{NSCH} for which the left-hand side of~\EQ{WeakTrac} is
defined, but the right-hand side is defined for a larger
class of solutions to~\EQ{NSCH} with weaker regularity; 
see also~\cite{Brummelen:2012fk,MelboKvamsdal2003,Zee:2011uq}.
The contribution of the fluid-solid surface tension in the right-member of~\EQ{dyn}
can be reformulated as:
\begin{equation}
\label{eq:surftens2}
-\int_{\Gamma_t}x\cdot\big(\sigma_{\mathrm{W}}(\varphi)\kappa{}n+\grad_{\Gamma}\sigma_{\mathrm{W}}(\varphi)\big)
=
\mathscr{C}(\hat{\chi}^{\mathrm{s}},\varphi;x)
+
\int_{\partial\Gamma_t}
\upsilon\cdot
\big(\sigma_{\mathrm{w}}(\varphi)\,x\big)
\end{equation} 
with 
\begin{equation}
\label{eq:Tdef}
\mathscr{C}(\hat{\chi}^{\mathrm{s}},\varphi;x)
=
\int_{\Gamma_t}
\sigma_{\mathrm{w}}(\varphi)\,
\grad_{\Gamma}\mathrm{Id}_{\Gamma_t}:
\grad_{\Gamma}x
\end{equation}
and $\mathrm{Id}_{\Gamma_t}$ the identity
on $\Gamma_t$ and $\upsilon$ the exterior unit normal vector to $\partial\Gamma_t$ in the
tangent bundle of~$\Gamma_t$; see~\cite{Bansch:2001fk,Dziuk:1991rz}. Let us note that
the right member of~\EQ{Tdef} depends implicitly on the solid deformation $\hat{\chi}^{\mathrm{s}}$ via the shape of the interface~$\Gamma_t$. The second term in the
right member of~\EQ{surftens2} cannot generally be bounded in weak formulations,
and it must vanish by virtue of boundary conditions on~$\hat{\chi}^{\mathrm{s}}$ or~$x$.

To provide a setting for the weak formulation of the 
fluid-solid-interaction problem, let $\smash[tb]{L^2(\omega)}$
denote the class of square-integrable functions on any $\omega\subset\IR^d$,
$H^1(\omega)$ the Sobolev space of functions in $L^2(\omega)$
with weak derivatives in $L^2(\omega)$, and $\smash[tb]{H^1_{0,\Xi}(\omega)}$
the subspace of functions that vanish on $\Xi\subseteq\partial\Omega$.
For a vector space~$X(\omega)$ of scalar-valued functions,
$\smash[tb]{X(\omega,\IR^d)}$ is the extension to the 
corresponding vector space of $\IR^d$-valued functions.
Given a vector space $V$ and a time interval $(0,T)$, 
$W(0,T;V)$ represents a (suitable) class of functions 
from $(0,T)$ into $X(\omega)$. 
%

We collect the ambient spaces for the fluid and solid variables 
into\footnote{The admissible solid deformations must in fact satisfy auxiliary conditions at the interface to ensure that 
the surface-tension contributions are well-defined. Detailed treatment of this aspect is beyond the scope of this work.}:
\begin{equation}
V:=
H^1_{0,\Gamma_{\mathrm{D}}^u}(\domF,\IR^d)
\times
{}L^2(\domF)
\times
H^1_{0,\Gamma_{\mathrm{D}}^{\varphi}}(\domF)
\times
H^1_{0,\Gamma_{\mathrm{D}}^{\mu}}(\domF)
\times
H^1_{0,\Gamma_{\mathrm{D}}^{\hat{\chi}}}(\domSref,\IR^d)
\end{equation}
For conciseness, we assume $g_{\mathrm{D}}^u|_{\Gamma_{\mathrm{D}}^u\setminus\Gamma_t}=0$, 
$g_{\mathrm{D}}^{\varphi}=0$, $g_{\mathrm{D}}^{\mu}=0$, and 
$\smash[tb]{g_{\mathrm{D}}^{\hat{\chi}}}=\mathrm{Id}$. 
The aggregated
fluid-solid-interaction problem can then be condensed into:
\begin{multline}
\label{eq:WeakFSI}
\text{{\em Find} }(u,p,\varphi,\mu,\hat{\chi}^{\mathrm{s}}-\mathrm{Id})\in{}
W(0,T;V)\text{ {\em such that almost everywhere in} }(0,T):
\\
d_t\langle\rho{}(u+\ell_{w|_{\Gamma_t}}),v+\ell_{x|_{\Gamma_t}}\rangle
+
\mathcal{A}_{\mathrm{N}}\big(u+\ell_{w|_{\Gamma_t}},w,\varphi;v+\ell_{x|_{\Gamma_t}}\big)
+
\mathcal{B}(p,v+\ell_{x|_{\Gamma_t}})
\\
+
\mathcal{B}(q,u+\ell_{w|_{\Gamma_t}})
+
d_t\langle\varphi,z\rangle
+
\mathcal{A}_{\mathrm{C}}(u+\ell_{w|_{\Gamma_t}},w,\varphi,\mu;z)
+
\mathcal{A}_{\mathrm{P}}(\varphi;\mu,y)
\\
+
d_{t}^2\langle\hat{\rho}\hat{\chi}^{\mathrm{s}},\hat{x}\rangle
+
\mathcal{W}'(\hat{\chi}^{\mathrm{s}},\hat{x})
+\mathscr{C}(\hat{\chi}^{\mathrm{s}},\varphi;x)
=
\mathcal{L}_{\mathrm{A}}(v,z,y,\hat{x})
\quad\forall(v,q,z,y,\hat{x})\in{}V 
\end{multline}
with $w=\partial_th_{\hat{\chi}^{\mathrm{s}}|_{\hat{\Gamma}}}\circ{}\hat{\chi}^{-1}$
and
$x=\hat{x}\circ{}\hat{\chi}^{-1}$, 
 and the 
aggregated linear form:
\begin{equation*}
\mathcal{L}_{\mathrm{A}}(v,z,y,\hat{x})
=
\int_{\Gamma_{\mathrm{N}}^u}v\cdot{}g_{\mathrm{N}}^u
+
\int_{\Gamma_{\mathrm{N}}^\mu}z\,g_{\mathrm{N}}^\mu
+
\int_{\Gamma_{\mathrm{N}}^\varphi}y\,g_{\mathrm{N}}^{\varphi}
+
\int_{\Gamma_{\mathrm{N}}^{\hat{\chi}}}\hat{x}\cdot{}g_{\mathrm{N}}^{\hat{\chi}}
\end{equation*}
It is to be noted that $\hat{\chi}^{\mathrm{s}}-\mathrm{Id}$ represents solid displacement.
Furthermore, by virtue of $(v+\ell_{x|_{\Gamma_t}})|_{\Gamma_t}=\hat{x}|_{\Gammaref}\circ{}\hat{\chi}^{-1}$ and
$(u+\ell_{w|_{\Gamma_t}})|_{\Gamma_t}=\partial_t\hat{\chi}|_{\Gammaref}\circ{}\hat{\chi}^{-1}$,
the test spaces for the equations of motion of the fluid and the solid and the 
trial spaces for the fluid and solid velocity in~\EQ{WeakFSI} are essentially continuous across the interface.
\vspace{-0.3\baselineskip}

\section{Numerical Experiments}
\label{sec:numexp}
\begin{figure}[t]
\includegraphics[width=\textwidth]{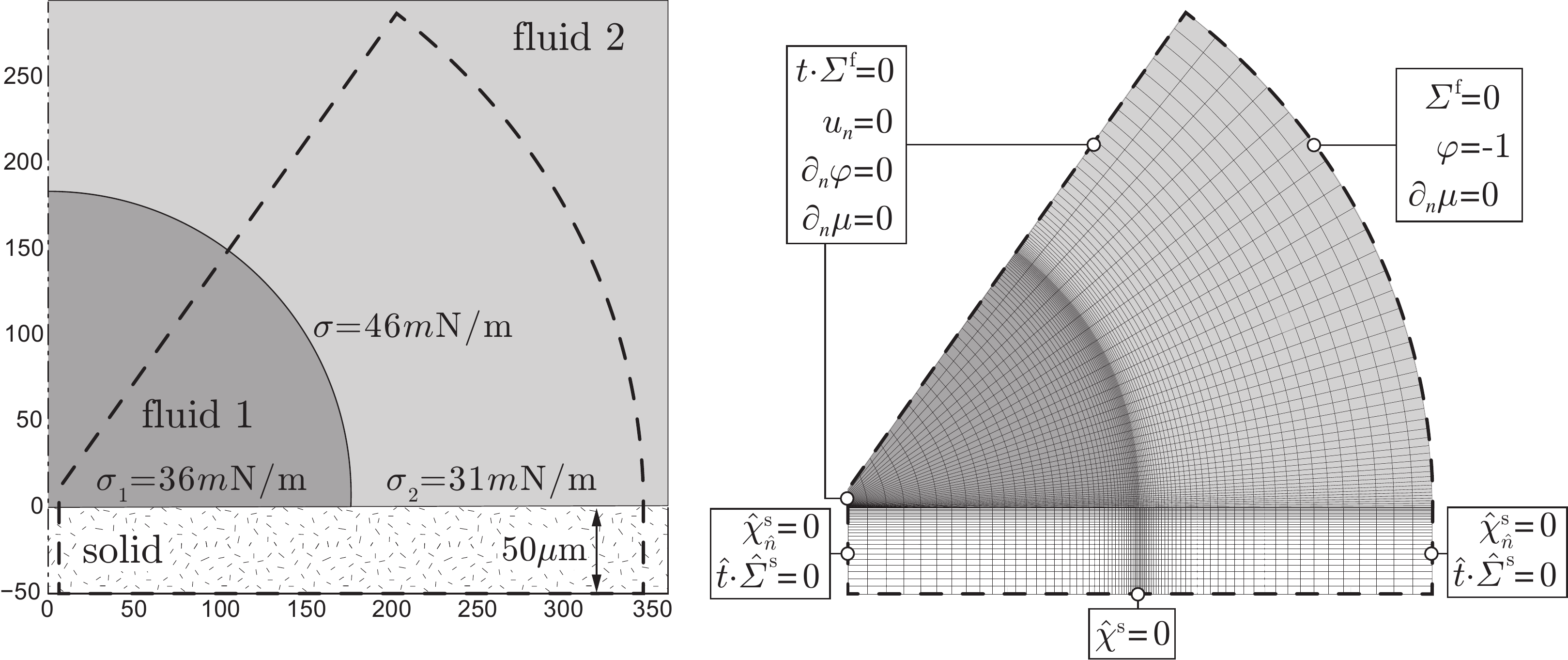}
\caption{Illustration of the considered experimental configuration ({\em left\/})
and the corresponding computational setup ({\em right\/}).\vspace{\baselineskip}
\label{fig:SetUp}}
\end{figure}

To evaluate the predictive capabilities of the presented CFSI model, we consider numerical approximations of~\EQ{WeakFSI} for the experimental setup in~\cite{Style:2013ay}. The test case  
concerns a $13.8\,p\mathrm{l}$ droplet on a soft substrate; see Fig.~\FIG{SetUp} ({\em left\/}). We characterize the substrate by a nearly incompressible solid with Saint Venant--Kirchhoff constitutive behavior, with Lam\'e parameters $\mu_{\mathrm{L}}=\tilde{E}/(2+2\tilde{\nu})$ and $\lambda_{\mathrm{L}}=\tilde{\nu}\tilde{E}/(1+\tilde{\nu})(1-2\tilde{\nu})$, and Young's modulus $\tilde{E}=3\,k\,\mathrm{Pa}$ and Poisson ratio $\tilde{\nu}=0.499$. 
The surface tension of the interface between the droplet (fluid 1) and ambient fluid (fluid 2) is $\sigma=46\,\mathit{m}\,\mathrm{N}/\mathrm{m}$. The fluid/solid surface tension of fluid 1 (resp. fluid 2) is $\sigma_1=36\,m\,\mathrm{N}/\mathrm{m}$ 
(resp. $\sigma_2=31\,m\,\mathrm{N}/\mathrm{m}$).
The diffuse-interface thickness is set to $\epsilon=2\,\mu\mathrm{m}$. Our interest is restricted to steady solutions and, hence, $\rho_1$, $\rho_2$, $\hat{\rho}$, $\nu$ 
and~$\gamma$ are essentially irrelevant.
For completeness, we mention that we select matched fluid 
densities $\rho=\rho_1=\rho_2=1.26\,p\mathrm{g}/(\mu\mathrm{m})^3$, 
matched fluid viscosities $\nu=\nu_1=\nu_2=1412\,m\,\mathrm{Pa}\,\mathrm{s}$,
solid density $\hat{\rho}=12.6\,p\mathrm{g}/(\mu\mathrm{m})^3$, and mobility
$\gamma=0.01\,(\mu\mathrm{m})^3\,\mu\mathrm{s}/p\mathrm{g}$.
We refer to Fig.~\FIG{SetUp} for further details of the experimental configuration.

We incorporate the rotational symmetry of the experimental setup 
in the discrete approximation of~\EQ{WeakFSI}. The considered approximation 
spaces are based on a locally refined mesh, adapted to the diffuse interface; 
see Fig.~\FIG{SetUp} ({\em right\/}). 
We apply Raviart-Thomas compatible B-spline approximations for 
velocity ($u$) and pressure~($p$) with
order $((3,2),(2,3))$ and $2$, respectively; see~\cite{Buffa2011,Evans:2012kx}. 
The order parameter ($\varphi$) and chemical potential ($\mu$)
 are approximated by means of quadratic B\nobreakdash-splines. 
 The solid deformation ($\hat{\chi}^{\mathrm{s}}$)
and the deformation of the fluid domain ($\smash[tb]{\hat{\chi}^{\mathrm{f}}}$) are approximated with quadratic B-splines as well. Let us note that
by virtue of the $C^1$\nobreakdash-continuity of the solid deformation, the
interface $\Gamma_t$ corresponds to a $C^1$ manifold. The temporal discretization
of~\EQ{WeakFSI} is based on backward Euler approximation of the time derivatives,
with time step $0.5\,m\mathrm{s}$.
In each time step, the aggregated fluid-solid interaction problem is solved
by means of subiteration with underrelaxation; see, for instance,~\cite{Brummelen:2011la}.

\begin{figure}[t]
\includegraphics[width=\textwidth]{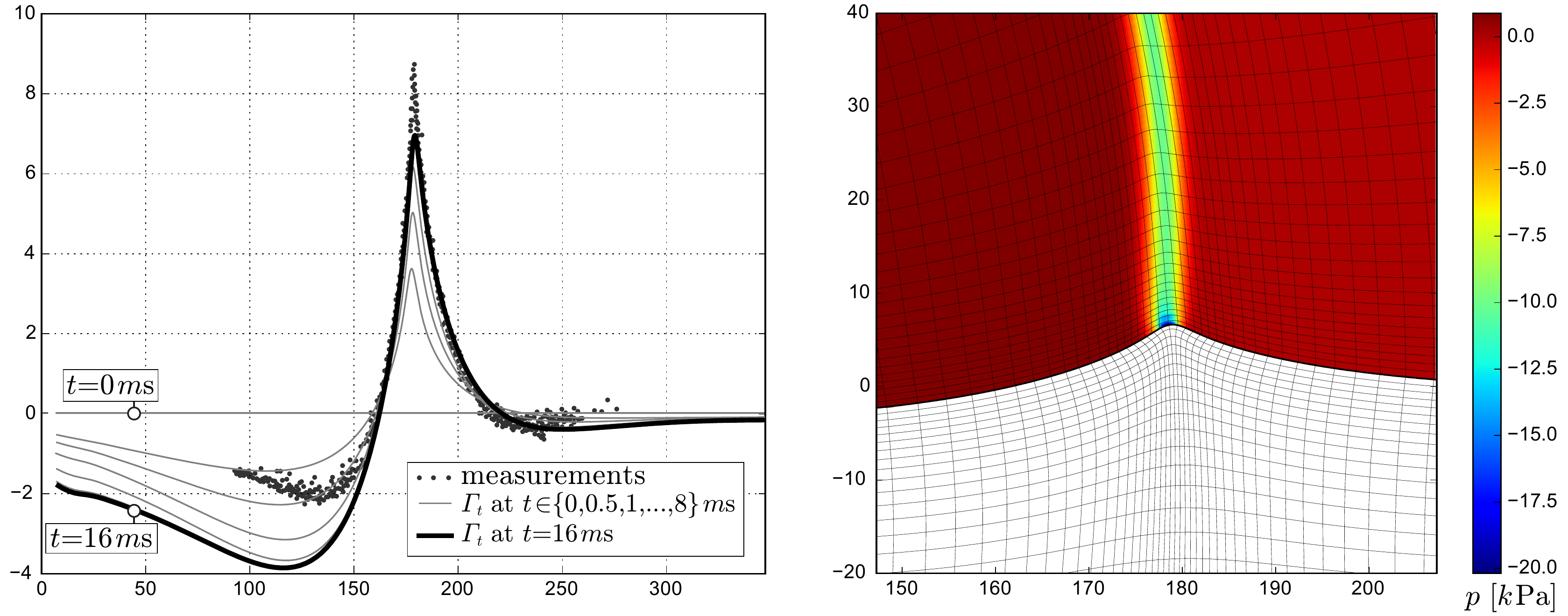}
\caption{Comparison of the computed fluid-solid-interface configuration $\Gamma_t$ at
$t\in\{0,0.5,1,2,4,8\}\,m\mathrm{s}$ ({\em grey\/}) and at $t=16\,m\mathrm{s}$ ({\em black\/}) and rendering of experimental results 
from~\cite{Style:2013ay} ({\em left\/}), and magnification of the contact-line region at~$t=16\,m\mathrm{s}$ with deformed fluid and solid meshes and computed pressure distribution~({\em right\/}).
\label{fig:Comparison}}
\end{figure}

Figure~\FIG{Comparison} ({\em left\/}) presents a comparison of the computed interface configuration, $\Gamma_t$,
at $t\in\{0,0.5,1,2,4,8,16\}\,m\mathrm{s}$ and experimental data from~\cite{Style:2013ay}. At $t=16\,m\mathrm{s}$, the interface has essentially reached its equilibrium deformation. The surface tension of the fluid-fluid interface yields a localized load on the fluid-solid interface near the contact line, resulting in a kink in the surface deformation of the soft substrate. In addition, the fluid-fluid surface tension 
leads to an increased pressure in the droplet relative to
the ambient pressure (see also Fig.~\FIG{Comparison} ({\em right\/})), viz. Laplace pressure, and a corresponding depression of the substrate. Comparison of the experimental and computed results conveys that the fluid-solid interface elevation at the contact line is underestimated by approximately 25$\%$. The underestimation can be attributed to the regularizing effect of the diffuse interface. It is anticipated that further reduction 
of the diffuse-interface thickness ($\epsilon$) and corresponding refinement of the mesh leads to an increase in the fluid-solid interface elevation at the contact line. 
The indentation of the substrate below the droplet is noticeably overestimated. In this regard, it is to be mentioned that
on account of the nearly incompressible behavior of the solid, its volume at $t=16\,m\mathrm{s}$ has decreased by only $0.16\%$ relative
to the initial volume.

Figure~\FIG{Comparison} ({\em right\/}) presents a magnification
of the contact-line region at $t=16\,m\mathrm{s}$ with the fluid and solid 
meshes in the actual configuration and the computed pressure distribution in 
the complex fluid. It is noteworthy that the pressure in the diffuse interface 
exhibits a localized minimum at the contact line. The pressure in the droplet is
virtually uniform with value $p\approx{}520\,\mathrm{Pa}$, which is close to the 
theoretical Laplace pressure $2\sigma/R$ in a droplet on a rigid substrate with radius $R=178\,\mu\mathrm{m}$ and surface tension $\sigma=46\,m\,\mathrm{N}/\mathrm{m}$.

\section{Conclusion}
\label{sec:concl}
We presented a model for the interaction of a complex fluid with an elastic solid, 
in which the complex fluid is represented by the Navier--Stokes--Cahn--Hilliard (NSCH)
equations and the solid is characterized by a hyperelastic material with a 
Saint Venant--Kirchhoff stored-energy functional.
The interaction between the fluid and the solid at their mutual interface is described 
by a preferential-wetting condition in addition to the usual kinematic and dynamic
interface conditions. The fluid traction on the fluid-solid interface comprises a non-standard capillary-stress contribution, and the dynamic condition contains
a contribution from the non-uniform fluid-solid surface tension. A weak formulation
of the complex-fluid-solid-interaction (CFSI) problem was presented, based on an ALE formulation of the 
NSCH system and a suitable reformulation of the complex-fluid traction and 
the fluid-solid surface-tension traction.

Numerical results were presented for a stationary droplet on a soft solid substrate, based on finite-element approximation of the weak formulation of the aggregated CFSI problem. Comparison of the computed results with experimental data for the
considered test case exhibited very good agreement in the contact-line region.
The substrate depression below the droplet was noticeably overestimated relative
to the experimental data. In view of the close agreement between the computed pressure
in the droplet and the theoretical Laplace pressure, it appears that the discrepancy
between the computed and observed depression is to be attributed to 
corresponding differences in the constitutive behavior of the solid substrate.
The overall good agreement between the computed and experimental data indicates
the potential of computational CFSI models based on the NSCH equations to 
predict elasto-capillary phenomena. 

\bibliography{BibFile}
\bibliographystyle{plain}

\end{document}